\documentclass[twoside]{article}
\usepackage{color}
\usepackage{amssymb}
\usepackage{dsfont}

\newtheorem{thm}{Theorem}[section]
\newtheorem{lem}[thm]{Lemma}

\newtheorem{pro}[thm]{Proposition}
\newtheorem{cor}[thm]{Corollary}
\newtheorem{axiom}{Definition}[section]
\newtheorem{rem}{Remark}[section]
\newtheorem{ex}{Example}[section]
\newtheorem{Ass}{Assumption}[section]
\newtheorem{notation}{Notation}[section]

\newcommand{\be}{\begin{equation}}
\newcommand{\ee}{\end{equation}}
\newcommand{\bde}{\begin{displaymath}}
\newcommand{\ede}{\end{displaymath}}
\newcommand{\beq}{\begin{eqnarray*}}
\newcommand{\eeq}{\end{eqnarray*}}
\newcommand{\beqa}{\begin{eqnarray}}
\newcommand{\eeqa}{\end{eqnarray}}
\newcommand{\bel }{\left\{\begin{array}{ll}}
\newcommand{\eel}{\cr \end{array} \right.}
\newcommand{\bd}{\begin{axiom} \rm }
\newcommand{\ed}{\end{axiom} \rm }
\newcommand{\brem }{\begin{rem} \rm }
\newcommand{\erem }{\end{rem}}
\newcommand{\bex}{\begin{ex} \rm }
\newcommand{\eex}{\end{ex}}
\newcommand{\bt}{\begin{thm}}
\newcommand{\et}{\end{thm}}
\newcommand{\bl}{\begin{lem}}
\newcommand{\el}{\end{lem}}
\newcommand{\bp}{\begin{pro}}
\newcommand{\ep}{\end{pro}}
\newcommand{\bcor}{\begin{cor}}
\newcommand{\ecor}{\end{cor}}

\newcommand{\droit}[1]{[#1]}
\newcommand{\oblique}[1]{\langle #1\rangle}
\newcommand{\abs}[1]{\left\vert#1\right\vert}

\def\proof{\noindent {\it Proof. $\, $}}
\def\finproof{\hfill $\Box$ \vskip 5 pt}

 \def\ind{\mathds{1}}

\newcommand{\bE}{\mathbb{E}}
\newcommand{\bF}{\mathbb{F}}
\newcommand{\bG}{\mathbb{G}}
\newcommand{\bH}{\mathbb{H}}

\newcommand{\bP}{\mathbb{P}}
\newcommand{\bQ}{\mathbb{Q}}
\newcommand{\bR}{\mathbb{R}}
\newcommand{\cA}{\mathcal{A}}
\newcommand{\cB}{\mathcal{B}}
\newcommand{\cE}{\mathcal{E}}
\newcommand{\cF}{\mathcal{F}}
\newcommand{\cG}{\mathcal{G}}
\newcommand{\cH}{\mathcal{H}}

\newcommand{\cL}{\mathcal{L}}
\newcommand{\cO}{\mathcal{O}}
\newcommand{\cP}{\mathcal{P}}

\newcommand{\V}{\vee}
\newcommand{\w}{\wedge}

\newcommand{\given}{\mathcal{\;\!|\;\!}}
\newcommand{\bgiven}{\mathcal{\;\!\big |\;\!}}

\def\hat{\widehat}
\def\tilde{\widetilde}

\begin{document}
\title{{\Large \bf Carthaginian Enlargement of Filtrations}
\footnote{This research was supported by the ``Chaire Risque de
Cr\'{e}dit'' of the French Banking Federation. All the authors are
members of the Laboratoire Analyse et Probabilit\'es of the
Universit\'e d'\'{E}vry-Val-D'Essonne, France}. \vskip 70 pt}

\author{
G. Callegaro\footnote{ Scuola Normale Superiore Pisa, Italy and
CREST, France. Email: giogiocallegaro@gmail.com} , M.
Jeanblanc\thanks{ Institut Europlace de Finance, France.
Email: monique.jeanblanc@univ-evry.fr} , B. Zargari\thanks{ Sharif
University of Technology, Iran. Email:
behnaz\_zargari@yahoo.com} }
\date{\today}
\maketitle
\begin{abstract}

This work is concerned with the theory of initial and progressive enlargements
of a reference filtration $\mathbb F$ with a random time $\tau$.
We provide, under an equivalence assumption, slightly stronger than the
absolute continuity assumption of Jacod, alternative proofs to results concerning canonical
decomposition of an $\bF$-martingale in the enlarged filtrations.
Also, we address martingales' characterization in the enlarged filtrations in
terms of martingales in the reference filtration, as well as predictable
representation theorems in the enlarged filtrations.
\end{abstract}



\section{Introduction} \label{4:Prelim}
We consider the case where a filtration $\mathbb F$ is enlarged
to give a filtration $\widetilde {\mathbb {F}}$, by means of a finite positive random variable $\tau$. In the literature,
two ways to realize such an enlargement are presented: either all of
a sudden at time $0$ (\emph{initial enlargement}), or
progressively, by considering the smallest filtration containing $\mathbb F$, satisfying the usual conditions, that makes $\tau$ a stopping time
(\emph{progressive enlargement}).

The ``pioneers'' who started exploring this research field,  at
the end of the seventies, were Barlow (see \cite{Ba78}), Jacod,
Jeulin and Yor (see the references that follow in the text).
The main questions that raised were the following: ``Does any $\mathbb
F$-martingale $X$ remain an $\tilde{\mathbb F}$-semimartingale?''
And, if it does: ``What is the semimartingale decomposition in
$\tilde{\mathbb F}$ of the $\mathbb F$-martingale $X$?''

The main contribution of the present work is to show how, under a
specific equivalence  assumption (see Assumption \ref{a:E_hypothesis}), slightly  stronger  than Jacod's
one in \cite{Jacod85}, some well-known fundamental results can be
proved in an alternative (and, in some cases, simpler) way. We
make precise that the goal of this paper is neither to present the
results in the most general case, nor to study the needed and
difficult regularity properties, for which we refer to existing
papers (e.g., Jacod \cite{Jacod85}).

Let us, now, motivate the title.
Inspired by a visit to the Tunisian archaeological site of
Carthage, where one can find remains of THREE levels of different
 civilizations, we decided to use the catchy adjective
``Carthaginian'' associated with filtration, since in this paper
there will be THREE levels of filtrations.

We consider, then,  three nested filtrations
$$
\bF \subset \bG \subset \bG^{\tau},
$$
where $\bG$ and $\bG^\tau$ stand, respectively, for the
\textit{progressive} and the \textit{initial} enlargement of $\bF$
with a finite random time  $\tau$ (i.e., a finite non-negative
random variable).

Under a specific assumption (see the $(\mathcal E)$-Hypothesis below), we address the
following problems:\vspace{-1mm}
\begin{itemize}
\item Characterization of $\bG$-martingales and $\bG^\tau$-martingales in terms of $\bF$-martingales;\vspace{-1mm}
\item Canonical decomposition of an $\bF$-martingale, as a semimartingale, in $\bG$ and $\bG^\tau$;\vspace{-1mm}
\item Predictable Representation Theorem in $\bG$ and $\bG^\tau$.
\end{itemize}

The exploited idea is the following: assuming that the
$\bF$-conditional law of $\tau$ is equivalent to the law of
$\tau$, after an \textit{ad hoc} change of probability measure,
the problem reduces to the   case where $\tau$ and $\mathbb F$ are
independent. Under this newly introduced probability measure,
working in the initially enlarged filtration  is ``easy''. Then,
under the original probability measure, for the initially enlarged
filtration, the results are achieved by means of Girsanov's
theorem. As for the progressively enlarged filtration, one
can proceed either by projecting on $\bG$ the results already
obtained for $\bG^{\tau}$ (e.g., in Proposition
\ref{p:pro_CharMart}), or, directly, by a change of probability
measure in the filtration $\bG$ (e.g., in Proposition \ref{p:PRP} (ii)).

The ``change of probability measure'' viewpoint for treating
problems on enlargement of filtrations was remarked in the early
80's and developed by Song in \cite{ShiqiTh}, and then by Ankirchner \emph{et al.} \cite{ADI}
(see also \cite{Jacod85}, Section 5). This is also  the  point of view
adopted by  Gasbarra \emph{et al.} in \cite{GVV2006} where the
authors apply the Bayesian approach to study the impact of the
initial enlargement of filtration on the characteristic triplet of
a semimartingale.

The paper is organized as follows. Section \ref{sect:Pre} introduces definitions and preliminary
results which will be crucial in the rest of the paper. Section \ref{Sect:martCaract} addresses the characterization
of $\bG$-martingales and $\bG^\tau$-martingales in terms of $\bF$-martingales. In Section \ref{Sect:canonDecomp},
the invariance of semimartingale property under the progressive and initial enlargements of filtration is studied,
and the formulae for the canonical decomposition of an $\bF$-martingale as a semimartingale in $\bG$ and $\bG^\tau$
are provided. In Section \ref{Sect:PRP}, we show that the enlarged filtrations $\bG$ and $\bG^\tau$ admit a
predictable representation property, as soon as the reference filtration $\bF$ enjoys one. Finally, Section
\ref{Sect:conclusion} ends the paper with some concluding remarks.

\section{Preliminaries}\label{sect:Pre}
We consider a probability space $(\Omega,\cA,\bP)$ equipped with a
filtration $\mathbb F={(\mathcal F_t)}_{t\geq 0}$ satisfying the
usual hypotheses of right-continuity and completeness, and where
$\mathcal F_0$ is the trivial $\sigma$-field, completed by the
$\bP$-negligible sets of $\cA$.

Let $\tau$ be a finite random time with law $\nu$,
$\nu(du)=\mathbb P(\tau\in du)$. We assume that $\nu$ has no atoms
and has $\bR^+$ as support.

We denote by $\cP(\bF)$ (resp. $\cO(\bF)$)
the predictable (resp. optional) $\sigma$-algebra corresponding to
$\bF$ on $\mathbb R^+ \times \Omega$.

Our standing assumption is the following:

\begin{Ass}\label{a:E_hypothesis}\textbf{$(\mathcal E)$-Hypothesis}\\
The $\mathbb F$-(regular) conditional law of $\tau$ is equivalent
to the law of $\tau$. Namely,
$$
\mathbb P(\tau \in du \vert \mathcal F_t) \sim \nu(du)\, \hspace{2mm} \textrm{for every} \ t \ge 0, \hspace{2mm} \mathbb P - a.s.
$$
\end{Ass}
This assumption, in the case where $t \in [0,T]$, corresponds to
the \textit{equivalence assumption} in F\"{o}llmer and Imkeller \cite{FoIm93} and in
Amendinger's thesis \cite[Assumption 0.2]{AmenTez}, and to
hypothesis (HJ) in the papers by Grorud and Pontier (see, e.g.,
\cite{GPInsider98}).

Amongst the consequences of the $(\mathcal E)$-Hypothesis,  one
has the existence and regularity of the conditional density, for
which we refer to Amendiger's reformulation (see remarks on page
17 of \cite{AmenTez}) of Jacod's result (Lemma 1.8 in
\cite{Jacod85}): there exists a strictly positive
$\cO(\bF)\otimes\cB(\bR^+)$-measurable function $(t, \omega, u)
\rightarrow p_t(\omega,u)$, such that for $\nu$-almost every $u \in \mathbb
R^+$, $p(u)$ is a c\`{a}dl\`{a}g $(\mathbb P,\mathbb
F)$-martingale  and
$$
\mathbb P(\tau >\theta \vert \mathcal F_t)= \int_\theta^\infty p_t(u)\nu (du)\quad  \textrm{for every} \; t \ge 0, \hspace{2mm} \mathbb P - a.s.
$$
In particular, $p_0(u)=1$ for $\nu$-almost every $u \in \mathbb R^+$. This family of processes $p$ is called the
 $(\mathbb P,\mathbb F)$-\emph{conditional density} of $\tau$ with
respect to $\nu$, or the \emph{density} of $\tau$ if there is no
ambiguity.

Furthermore, under the $(\mathcal E)$-Hypothesis, the assumption
that $\nu$ has no atoms implies that the default time $\tau$
avoids the $\mathbb F$-stopping times, i.e., $\mathbb P(\tau
=\xi)=0$ for every $\mathbb F$-stopping time $\xi$ (see, e.g.,
Corollary 2.2 in El Karoui \emph{et al.} \cite{ejj}).

The initial enlargement of $\bF$ with $\tau$, denoted by $\mathbb
G^{\tau}=(\mathcal G^\tau_t, t\geq 0)$, is defined as
$\cG^\tau_t=\cF_{t }\V \sigma(\tau)$. It was shown in \cite[Proposition 1.10]{AmenTez}  that the strict positiveness of
$p$ implies the right-continuity of the filtration $\bG^\tau$.

Let $\bH={(\cH_t)}_{t\geq 0}$ denote the smallest filtration with
respect to which $\tau$ is a stopping time, i.e.,
$\cH_t=\sigma(\ind_{\tau \leq s},s\leq t)$. This filtration is
right-continuous. The progressive enlargement of $\bF$ with the random time
$\tau$, denoted by $\bG={(\cG_t)}_{t\geq 0}$, is defined as the
right-continuous regularization of  $ \bF \vee \bH$.

In the sequel, we will consider the right-continuous version of
all the martingales.

Now, we consider the change of probability measure introduced,
independently, by Grorud and Pontier in \cite{GPInsider98} and by
Amendinger in \cite{AmenTez}. Having verified that the process
$L$, given by $L_t=\frac{1}{p_t(\tau)}$, $t \ge 0$, is a
$(\bP,\bG^{\tau})$-martingale, with $\mathbb E(L_t)=L_0=1$, these
authors  defined a locally equivalent probability measure $\bP^*$
setting
$$
{d\bP^*}_ {\vert {\cG^{\tau}_t}}= L_t\; {d\bP}_{ \vert{\cG^{\tau}_t}} = \frac{1}{p_t(\tau)} {d\bP}_{ \vert{\cG^{\tau}_t}} \;.
$$
They proved that, under $\bP^*$, the random time  $\tau$ is
independent of $\mathcal F_t$ for any $t \ge 0$ and, moreover,
that
$$
{\bP^*}_{\vert{\cF_t}} = {\bP}_{\vert {\cF_t}}\;\; \textrm{for any}\ t \ge 0, \; \quad {\bP^*}_{\vert {\sigma(\tau)}}= {\bP}_{\vert {\sigma(\tau)}}.
$$
The above properties imply that $\mathbb P^* (\tau \in du \vert
\mathcal F_t) = \mathbb P^*(\tau \in du)$, so that the $(\mathbb
P^*,\mathbb F)$-density of $\tau$, denoted  by  $p^*(u), u \ge 0$,
is a constant equal to one, $\mathbb P^* \otimes \nu$-a.s.

\begin{rem}
The probability measure $\bP^*$, being defined on $\cF_t$ for $t\geq 0$,
is (uniquely) defined on $\cF_\infty= {\bigvee}_{t\geq 0}\cF_t$.
Then, as $\tau$ is independent of $\bF$ under $\bP^*$, it immediately follows that $\tau$
is also independent of $\cF_\infty$, under $\bP^*$. However, one can not claim that:
``$\mathbb P^*$ is equivalent to $\bP$ on $\cG^{\tau}_{\infty}$'', since we do not
know \emph{a priori} whether $\frac{1}{p(\tau)}$ is a \emph{closed} $(\bP,\bG^\tau)$-martingale or not.
A similar problem is studied by F\"{o}llmer
and Imkeller in \cite{FoIm93} (it is therein called ``paradox'')
in the case  where  the reference (canonical) filtration is
enlarged by means of the information about the endpoint at time $t=1$. In
our setting, it corresponds to the case where $\tau \in
\cF_\infty$ and $\tau \notin \cF_t, \forall \ t$. \end{rem}

\begin{rem}\label{r:mart}
Let $x=(x_t,t\geq 0)$ be a $(\bP,\bF)$-local martingale.
Since $\bP$ and $\bP^*$ coincide on $\bF$, $x$ is a
$(\bP^*,\bF)$-local martingale, hence, using the fact that
$\tau$ is independent of $\mathbb F$ under $\bP^*$, a
$(\bP^*,\bG)$-local martingale (and also a
$(\bP^*,\bG^\tau)$-local martingale).
\end{rem}

\begin{notation}
In this paper, as we mentioned, we deal with three different
levels of information and two equivalent probability measures. In
order to distinguish objects defined under $\bP$ and under
$\bP^*$, we will use a superscript $*$ when working under $\bP^*$.
For example, $\mathbb E$ and $\mathbb E^*$ stand for the
expectations under $\mathbb P$ and $\mathbb P^*$, respectively.
For what concerns the filtrations, when necessary, we will use the
following illustrating notation: $x, X, X^\tau$ to denote
processes adapted to $\bF, \bG$ and $\bG^{\tau}$, respectively
(we shall not use the same notation for processes stopped at $\tau$,
so that there will be no possible confusion for the notation $X^\tau$).
\end{notation}

The following proposition provides a relation between martingales w.r.t.
a ``larger'' and a ``smaller'' filtration and turns out to be useful in the subsequent sections:
\begin{pro}\label{p:project}
Let  $\tilde\bF$ be a filtration larger than $\bF$, that is
$\cF_t\subseteq\tilde\cF_t$, for every $t\geq 0$. If $x$ is a u.i. (uniformly integrable)
$\bF$-martingale, then there exists an $\tilde\bF$-martingale
$\tilde x$,   such that $\bE(\tilde x_t\given\cF_t)=x_t$, $t\geq
0$.
\end{pro}
\proof The process $\tilde x$ defined by $\tilde
x_t:=\bE(x_\infty\given\tilde\cF_t)$ is an $\tilde
\bF$-martingale, and
$$
\bE(\tilde x_t\given\cF_t)=\bE\big(\bE(x_\infty
\given\tilde\cF_t)\given\cF_t\big)= \bE(x_\infty \given\cF_t)
=x_t\;.
$$
\finproof
\begin{rem} The uniqueness of such a martingale $\tilde x$ is not
claimed in the above proposition and does not hold in general.
\end{rem}
We now recall the definition of the \emph{immersion property}, a terminology which will be used in the sequel.
 Under a given probability measure $\bQ$, a filtration $\bF$ is said to be \emph{immersed} in a larger filtration $\tilde\bF$, if every $(\bQ,\bF)$-martingale is a $(\bQ,\tilde\bF)$-martingale.

\subsection{Characterization of different measurability properties }
Before focusing on the three topics announced from the beginning,
we recall some important results on the characterization of
$\mathcal G_t^ \tau$ and $\mathcal G_t $-measurable random
variables, as well as  $\bG^{\tau}$ and $\bG$-predictable
processes. The necessary part of the result below, in the case of
predictable processes, is due to Jeulin    \cite[Lemma 3.13]{Jeulin80}. See also Yor \cite{YorNoyaux85}.
\begin{pro}\label{p:ini_categorize} One has
\begin{itemize}
\item [(i)] A random variable $Y_t^{\tau}$ is $\mathcal G_t^\tau$-\emph{measurable} if and only if it is of the form $Y_t^{\tau}(\omega)=y_t(\omega,\tau(\omega))$, for some
$\mathcal F_t \otimes \cB(\mathbb R^+)$-measurable random variable
$y_t(\cdot,u)$.
\item[(ii)] A process $Y^{\tau}$ is $\bG^\tau$-\emph{predictable} if and only if it is of
the form $Y_t^{\tau}(\omega)=y_t(\omega,\tau(\omega))$, $t \ge 0$, where
$(t,\omega,u)\mapsto y_t(\omega,u)$ is a $\cP(\bF)\otimes\cB(\mathbb R^+)$-measurable function.
\end{itemize}
\end{pro}
\proof The proof of part (i) is based on the fact that  $\mathcal
G_t^{\tau}$-measurable random variables are generated by random
variables of the form $X_t(\omega)=x_t(\omega) f \big(
\tau(\omega) \big)$, with $x_t \in \mathcal F_t$ and $f$ bounded
Borel function on $\mathbb R^+$.

(ii) It suffices to notice that  processes of the form $X_t :=x_t
f(\tau) $, $t \ge 0$, where $x $  is  $ \mathbb F$-predictable and
$f$  is a bounded Borel function on $\mathbb R^+$, generate the
$\mathbb G^{\tau}$-predictable $\sigma$-field.
 \finproof

For what concerns the progressive enlargement setting, the
following result is analogous to Proposition
\ref{p:ini_categorize}. The necessity of part (ii) is already
proved in Jeulin   \cite[Lemma 4.4 ]{Jeulin80}.
\begin{pro}\label{p:pro_categorize}One has
\begin{itemize}
\item [(i)] A random variable $Y_t$ is $\mathcal G_t$-measurable
if and only if it is of the form $Y_t(\omega)=\tilde y_t(\omega)
\ind_{t<\tau(\omega)}+\hat y_t(\omega,\tau(\omega))
\ind_{\tau(\omega) \leq t}$ for some $\mathcal F_t$-measurable
random variable $\tilde y_t$ and  some   family of  $\mathcal F_t
\otimes \cB(\mathbb R^+)$-measurable random variables  $\hat
y_t(\cdot,u),t\geq u$.
\item[(ii)] A process $Y$ is $\bG$-predictable if and only if it is of the form
$Y_t(\omega)=\tilde y_t(\omega) \ind_{t \le \tau(\omega)} + \hat
y_t(\omega,\tau(\omega))\ind_{\tau(\omega) < t}$, $t \ge 0$, where $\tilde y$
is $\mathbb F$-predictable and $(t,\omega,u)\mapsto \hat y_t(\omega,u)$ is a
$\cP(\bF)\otimes\cB(\mathbb R^+)$-measurable function.
\end{itemize}
\end{pro}
\proof For part (i), it suffices to recall that $\mathcal
G_t$-measurable random variables are generated by random variables
of the form $X_t(\omega)=x_t(\omega) f \big( t \wedge \tau(\omega)
\big)$, with $x_t \in \mathcal F_t$ and $f$ a bounded Borel
function on $\mathbb R^+$.

(ii) It suffices to notice that $\mathbb G$-predictable processes
are generated by processes of the form $X_t =x_t  \ind_{t\leq \tau
}+  \hat x_t f(\tau) \ind_{  \tau <t }$, $t \ge 0$, where   $x,
\hat x$  are $\mathbb F$-predictable and $f$ is a bounded Borel
function, defined on $\mathbb R^+$. \finproof

Such a characterization result does \emph{not} hold for optional
processes, in general. We refer to Barlow \cite[Remark on pages
318 and 319]{Ba78}, for a counterexample.

\subsection{Expectation and projection tools}

\begin{lem}\label{l:ini_expect}
Let $Y_t^{\tau}=y_t(\tau)$ be a $\mathcal G_t^{\tau}$-measurable random variable.\\
(i) If  $y_t(\tau)$ is $\bP$-integrable and  $y_t(\tau)=0\ \bP$-a.s. then, for $\nu$-a.e. $u \ge 0$, $y_t(u)=0$ $\bP$-a.s.\\
(ii) For $s \le t$ one has, $\bP$-a.s. (or, equivalently,
$\bP^*$-a.s.): \\ if $y_t(\tau)$ is $\bP^*$-integrable and if
$y_t(u)$  is $\bP$ (or $\bP^*$)-integrable for any $u \ge 0$,
\begin{equation} \label{11}
\mathbb E^* \big(y_t(\tau)\given\cG^{\tau}_s\big) = {\mathbb E^*
\big(y_t(u)\given \cF_s\big)}_{\big \vert u = \tau} =
 {\mathbb E \big(y_t(u)\given \cF_s\big)}_{\big \vert u =
 \tau};
 \end{equation}
 if $y_t(\tau)$ is $\bP$-integrable
\begin{equation}
\bE \big(y_t(\tau)\given\cG^{\tau}_s\big)=\frac{1}{p_s(\tau)}
{\bE\big(y_t(u)p_t(u)\given\cF_s\big)}_{\big \vert_{u=\tau}} \; .
\label{12}\end{equation}
\end{lem}
\proof (i) We have, by applying Fubini-Tonelli's Theorem,
\begin{eqnarray*}
0 & = & \bE\big(\abs{y_t(\tau)}\big)
 =\bE\Big(\bE\big(\abs{y_t(\tau)} \big\vert \cF_t\big)\Big) = \bE\Big(\int_0^\infty \abs{y_t(u)} p_t(u) \nu(du) \Big)  \,.
\end{eqnarray*}
  Then $\int_0^\infty \abs{y_t(u)} p_t(u)
\nu(du)=0$ $\mathbb P$-a.s. and, given that $p_t(u)$ is strictly
positive for any $u$, we have that, for $\nu$-almost every $u$, $y_t(\cdot,u)=0$ $\bP$-a.s.
\\
(ii) The  first equality in (\ref{11}) is straightforward for
elementary random variables of the form $f(\tau) x_t$, given the
independence between $\tau$ and  $\mathcal F_{t}$, for any $t \ge
0$. It is extended  to $\mathcal G_t^{\tau}$-measurable r.v's via
the monotone class theorem. The second equality follows from the
fact that $\mathbb P$ and $\mathbb P^*$ coincide on $\mathcal
F_t$, for any $t \ge 0$.

Equality (\ref{12}) is an immediate consequence  of
(\ref{11}), since it suffices, by means of (conditional) Bayes'
formula, to pass under the measure $\mathbb P^*$. Namely, for
$s<t$, we have
$$
\bE\big(y_t(\tau)\given\cG^\tau_s\big)=
\frac{\bE^*\big(y_t(\tau)p_t(\tau)\given\cG^\tau_s\big)}{\bE^*\big(p_t(\tau)\given\cG^\tau_s\big)}=
\frac{1}{p_s(\tau)} {\bE\big(y_t(u)p_t(u)\given\cF_s\big)}_{\bgiven_{u=\tau}},
$$
where in the last equality we have   used  the previous result
(\ref{11}) and the fact that $p(\tau)$ is a $(\mathbb P^*,\mathbb
G^{\tau})$-martingale.  Note that if $y_t(\tau)$ is $\mathbb
P$-integrable, then $\bE( \int_0^\infty  \vert y_t(u) \vert p_t(u)
\nu(du))= \bE(\vert y_t(\tau)\vert)<\infty$, which implies that
$\bE(\vert y_t(u) \vert p_t(u) )<\infty$. \finproof

When working with the progressively enlarged filtration $\bG$, it is convenient to introduce
the notation  $G$ (resp., $G^*(\cdot)$) for the Az\'ema supermartingale associated with $\tau$
under the probability measure $\bP$ (resp., $\bP^*$).
More precisely,
\begin{eqnarray}
&& G_t := \mathbb P( \tau > t \vert \cF_t)= \int_t^\infty p_t(u)\nu (du),\;\label{Eq:G} \\
&& G^*(t):= \mathbb P^*( \tau > t \vert \cF_t)= \bP^*(\tau >t)=\bP(\tau >t)=\int_t^\infty \nu (du). \label{Eq:G(t)}
\end{eqnarray}
Note, in particular,  that ${(G_t)}_{t\geq 0}$ is an $\bF$
supermartingale, whereas $G^*(\cdot)$ is a (deterministic)
continuous and decreasing function. Furthermore, it is clear that,
under the $(\mathcal E)$-Hypothesis and the  hypothesis that
the support of $\nu$ is $\bR^+$, $G$ and $G^*(\cdot)$ do not
vanish.
\begin{lem}\label{l:pro_expect}
Let $Y^{\tau}_t=y_t(\tau)$ be a $\cG_t^\tau$-measurable,
$\bP$-integrable random variable. Then, for $s \le t$,
$$
\bE(Y_t^{\tau}\given\cG_s)=\bE( y_t(\tau)\given\cG_s)=\tilde
y_s\ind_{s<\tau}+\hat y_s(\tau)\ind_{\tau\leq s},
$$
with
\begin{eqnarray*}
\tilde y_s &=& \frac{1}{G_s}\bE\Big(\int_s^{+\infty} y_t(u)p_t(u)\nu(du)\given\cF_s\Big)\;,\\
\hat y_s(u)&=& \frac{1}{p_s(u)}\bE\big(y_t(u)p_t(u)\given\cF_s\big)\;.
\end{eqnarray*}
\end{lem}
\proof From the above Proposition \ref{p:pro_categorize}, it is
clear that $\bE(y_t(\tau)\given\cG_s)$ can be written in the form
$\tilde y_s\ind_{s<\tau}+\hat y_s(\tau)\ind_{\tau\leq s}$. On the
set
 $\{s<\tau\}$, we have, applying Lemma 3.1.2 in Bielecki \emph{et al.} \cite{BJR} and using the $(\mathcal E)$-Hypothesis (see also \cite{ejj} for analogous computations),
\begin{eqnarray*}
    \ind_{s<\tau} \bE(y_t(\tau)\given\cG_s)
&=& \ind_{s<\tau} \frac{\bE\left[\bE(y_t(\tau) \ind_{s<\tau} \given \cF_t)\vert \cF_s \right]}{G_s}\\
&=& \ind_{s<\tau} \frac{1}{G_s} \bE \left( \int_s^{+ \infty} y_t(u) p_t(u) \nu(du) \bgiven\cF_s \right)
 =: \ind_{s<\tau} \tilde y_s.
\end{eqnarray*}
On the complementary set, we have, by applying Lemma
\ref{l:ini_expect},
$$
\ind_{\tau \le s} \bE(y_t(\tau)\given\cG_s) = \ind_{\tau \le s} \bE\left[\bE(y_t(\tau)\given \cG_s^{\tau})\vert \cG_s \right] = \ind_{\tau \le s} \frac{1}{p_s(\tau)} {\bE\big(y_t(u)p_t(u)\given\cF_s\big)}_{\big \vert_{u=\tau}}  =: \ind_{\tau \le s} \hat y_s(\tau).
$$
\finproof

For $s>t$, we obtain $\bE(Y_t^{\tau}\given\cG_s) = \frac 1{G_s}
\int_s^\infty y_t(u)p_s(u) \nu(du) \ind_{s<\tau}+
y_t(\tau)\ind_{\tau\leq s}$.

As an application,  projecting  the martingale $L$ (defined
earlier as $L_t=\frac{1}{p_t(\tau)}, t \ge 0$) on $\bG$ yields to
the corresponding Radon-Nikod\'ym density on $\mathbb G$:
$$
{d\bP^*}_{ \vert {\cG_t}}= \ell_t \ {d\bP}_{ \vert{\cG_t}}\;,
$$
with
\begin{eqnarray*}
\ell_t &:=& \bE(L_t \vert \cG_t)   =
\ind_{t<\tau}\frac{1}{G_t}\int_t^\infty \nu(du) + \ind_{\tau
\leq t}\frac{1}{p_t(\tau)}\\&=&
\ind_{t<\tau}\frac{G^*(t)}{G_t}+ \ind_{\tau \leq t
}\frac{1}{p_t(\tau)}\;.
\end{eqnarray*}

We now recall some useful facts concerning the compensated
martingale of the default indicator process $H_t=\ind_{ \tau
\leq t }, t \ge 0$. We know, from the general theory (see, for
example, \cite{ejj}), that the process $M$ defined as
\begin{equation}\label{Eq:M}
M_t:=H_t-\int_ 0^{t\wedge \tau} \lambda _s \ \nu(ds), \quad t \ge 0,
\end{equation}
with $\lambda _t =\frac{p_t(t)}{G_t }$, is a $(\bP,\bG)$-martingale and that
\begin{equation}\label{Eq:M*}
M^*_t := H_t - \int_0^{t\wedge \tau} \lambda^*(s) \ \nu(ds), \quad t \ge 0,
\end{equation}
with $\lambda^*(t)=\frac{1}{G^*(t)}$, is a $(\bP^*,\bG)$-martingale.
Furthermore,  since  $\lambda^*$  is  deterministic, $M^*$ (being
$\mathbb H$-adapted) is a $(\bP^*,\bH)$-martingale, too.

 We conclude this subsection with the following two
propositions, concerning the predictable projection, respectively
on $\mathbb F$ and on $\mathbb G$, of a $\mathbb
G^{\tau}$-predictable process. The first result is due to Jacod
 \cite[Lemma 1.10]{Jacod85}.
\begin{pro}\label{p:predictable_proF}
Let $Y^\tau=y(\tau)$ be a $\bG^\tau$-predictable, positive or
bounded, process. Then, the $\bP$-predictable projection of
$Y^\tau$ on $\mathbb F$ is given by
$$
^{(p)}{(Y^{\tau})}_t = \int_0^\infty y_{t}(u) p_{t-}(u) \nu(du) \;.
$$
\end{pro}
\proof It is obtained by a monotone class argument and by using
the definition of density of $\tau$, writing, for ``elementary''
processes, $Y_t^{\tau} : = y_t  f(\tau)$, with $y$ a bounded
$\mathbb F$-predictable process and $f$ a bounded Borel function.
For this, we refer to the proof of Lemma 1.10 in \cite{Jacod85}.
\finproof
\begin{pro}\label{p:predictable_pro}
Let $Y^\tau=y(\tau)$ be a $\bG^\tau$-predictable, positive or
bounded, process. Then, the $\bP$-predictable projection of
$Y^\tau$ on $\bG$ is given by
$$
^{(p)}{(Y^{\tau})}_t=\ind_{t\leq\tau}\frac{1}{G_{t-}}\int_t^\infty y_{t}(u)p_{t-}(u)\nu(du)+\ind_{\tau< t}y_{t}(\tau)\;.
$$
\end{pro}
\proof In this proof, for clarity, the left-hand side superscript
``$(p \, \mathbb G )$'' denotes the $\bP$-predictable projection
on $\mathbb G$, while the left-hand side superscript ``$(p \,
{\mathbb F} )$'' indicates the $\bP$-predictable projection on
$\mathbb F$. By the definition of predictable projection, we know
(from Proposition \ref{p:pro_categorize} (ii)) that we are looking
for a (unique) process of the form
$$
^{(p\, {\mathbb G})}{(Y^{\tau})}_t = \tilde y_t \ind_{t \le \tau}
+ \hat y_t(\tau)\ind_{\tau < t}, \quad t \ge 0,
$$
where $\tilde y$ is $\mathbb F$-predictable, positive or bounded, and
$(t,\omega,u)\mapsto \hat y_t(\omega,u)$ is a $\cP(\bF)\otimes\cB(\mathbb
R^+)$-measurable positive or bounded function, to be identified.
\begin{itemize}
\item On the predictable set $\{ \tau < t \}$, being $Y^{\tau}$ a $\mathbb G^{\tau}$-predictable, positive or bounded, process (recall Proposition \ref{p:ini_categorize} (ii)), we immediately find $\hat y(\tau) = y(\tau)$;
\item On the complementary set $\{ t \le \tau \}$, introducing the $\mathbb G$-predictable process
$$
Y:=\,  ^{(p \,{\mathbb G})} \hspace{-0.7 mm}({Y^{\tau}})
$$
it is possible to use Remark 4.5, page 64 of Jeulin \cite{Jeulin80} (see
also Dellacherie and Meyer \cite[Ch. XX, page 186]{DellMeyXXVII-XXIV}), to write
$$
  Y\ind_{ ]\!] 0 , \tau ]\!]} = \frac{1}{G_{-}} \ \ ^{(p \,{\mathbb F})} {\left(Y  \ind_{ ]\!] 0 , \tau ]\!]} \right)}
    \ind_{ ]\!] 0 , \tau ]\!]} = \frac{1}{G_{-}} \ \ ^{(p\, {\mathbb F})} {\left( ^{(p\, {\mathbb G})} \hspace{-0.3 mm}{(Y^{\tau})} \ind_{ ]\!] 0 , \tau ]\!]} \right)}  \ind_{ ]\!] 0 , \tau ]\!]}.
$$
We then have, being $\ind_{ ]\!] 0 , \tau ]\!]}$, by definition, $\mathbb G$-predictable (recall that $\tau$ is a $\mathbb G$-stopping time),
$$
  Y\ind_{ ]\!] 0 , \tau ]\!]} = \frac{1}{G_{-}} \ \ ^{(p \,{\mathbb F})} {\left( {Y^{\tau} \ind_{ ]\!] 0 , \tau ]\!]}}  \right)}  \ind_{ ]\!] 0 , \tau ]\!]},
$$
where the last equality follows by the definition of predictable
projection, being $\mathbb F \subset \mathbb G$. Finally, given
the result in Proposition \ref{p:predictable_proF}, we have
$$
^{(p\,{\mathbb F})} {\left( {Y^{\tau} \ind_{ ]\!] 0 , \tau ]\!]}}
\right)}_t = \int_{t}^{+ \infty} y_{t}(u) p_{t-}(u) \nu(du)
$$
and the proposition is proved.
\end{itemize}
\finproof
\section{Martingales' characterization}\label{Sect:martCaract}
The aim of this section is to characterize $(\mathbb P,\mathbb
G^{\tau})$ and $(\mathbb P,\mathbb G)$-martingales in terms of
$(\mathbb P,\mathbb F)$-martingales.

\begin{pro}\label{p:ini_CharMart}
\textbf{Characterization of $(\bP,\bG^\tau)$-martingales in terms of $(\bP,\bF)$-martingales}\\
A process $Y^{\tau}=y(\tau)$ is a $(\bP,\bG^\tau)$-martingale if
and only if   $(y_t(u)p_t(u),t\geq 0)$ is a
$(\bP,\bF)$-martingale, for $\nu$-almost every $u \ge 0$.
\end{pro}
\proof
The sufficiency is a direct consequence of Proposition \ref{p:ini_categorize} and Lemma \ref{l:ini_expect} (ii).\\
Conversely, assume that $y(\tau)$ is a $\bG^\tau$-martingale. Then, for $s \le t$, from Lemma \ref{l:ini_expect} (ii),
$$
y_s(\tau)=\mathbb E \left( y_t(\tau) \vert \mathcal G_s^{\tau} \right) =
\frac{1}{p_s(\tau)} {\bE\big(y_t(u)p_t(u)\given\cF_s\big)}_{\vert {u=\tau}}
$$
and the result follows from Lemma \ref{l:ini_expect} (i).
\finproof
\begin{rem}\label{r:ini_CharMartLoc}
This result, being a consequence of the Girsanov theorem (cf. the proof of Lemma \ref{l:ini_expect} (ii)), can immediately be extended to $(\bP,\bG^\tau)$-local martingales.\end{rem}

Passing to the progressive enlargement setting, we state and prove
a martingale characterization result, first formulated by El
Karoui \emph{et al.} in \cite[Theorem 5.7]{ejj}.

\begin{pro}\label{p:pro_CharMart}
\textbf{Characterization of $(\bP,\bG)$-martingales in terms of $(\bP,\bF)$-martingales}\\
A  $\bG$-adapted process $Y_t:=\tilde y_t\ind_{t<\tau}+\hat y_t(\tau)\ind_{\tau\leq t}, t \ge 0$,
is a $({\bP},\bG)$-martingale if and only if the following two conditions are satisfied
\begin{itemize}
\item [(i)] for $\nu$-almost every $u\geq 0$, $\big(\hat y_t(u)p_t(u), t\geq
u\big)$ is a $(\bP,\bF)$-martingale; \item [(ii)] the process
$m=(m_t,t\geq 0)$, given by
\begin{equation}\label{e:m}
m_t:=\bE(Y_t\given\cF_t)=\tilde y_tG_t+\int_0^t\hat y_t(u)p_t(u)\nu(du)\;,
\end{equation}
is a $(\bP,\bF)$-martingale.
\end{itemize}
\end{pro}

\proof For the necessity, in a first step, we show that we can
reduce our attention to the case where $Y$ is u.i.: indeed, let
$Y$ be a $(\bP,\bG)$-martingale. For any $T$, let $Y^{(T)}
=(Y_{t\w T},t\geq 0)$  be  the associated stopped martingale,
which is u.i. Assuming that the result is established for u.i.
martingales will prove that  the processes in \textit{(i)} and
\textit{(ii)} are martingales up to time $T$. Since $T$ can be
chosen as large as possible, we shall have the result.

Assume,  then, that $Y$ is a u.i. $(\bP,\bG)$-martingale. From
Proposition \ref{p:project}, $Y _{t }=\bE(Y^{\tau }_{t
}\given\cG_t)$ for some $(\bP,\bG^{\tau})$-martingale $Y^{\tau }$.
Proposition \ref{p:ini_CharMart}, then, implies that $Y^{\tau }_{t
}= y_{t }(\tau)$, where for $\nu$-almost every $u \ge 0$ the
process $\big(y_{t}(u)p_t(u),t \ge 0\big)$ is a
$(\bP,\bF)$-martingale. One then has
$$
\ind_{\tau\leq t}\hat y_t(\tau)=\ind_{\tau\leq t}Y_t=
\ind_{\tau\leq t}\bE(Y^{\tau }_{t }\given\cG_t)=\bE(\ind_{\tau\leq t}Y^{\tau }_{t }\given\cG_t)=
\ind_{\tau\leq t}y_t(\tau)\;,
$$
which implies, in view of Lemma \ref{l:ini_expect}(i), that for
$\nu$-almost every $u\leq t$, the identity $y_t(u)= \hat y_t(u)$
holds $\bP$-almost surely. So, \textit{(i)} is proved.
 Moreover, $Y$ being a $(\bP,\bG)$-martingale, its projection
on the smaller filtration $\bF$, namely the process $m$ in
(\ref{e:m}), is a $(\bP,\bF)$-martingale.

Conversely, assuming \textit{(i)} and \textit{(ii)}, we verify that $\bE(Y_t\given\cG_s)=Y_s$ for $s \le t$.
We start by noting that
\begin{equation}\label{e:pro_expect}
\bE(Y_t\given\cG_s)=\ind_{s<\tau}\frac{1}{G_s}\bE(Y_t\ind_{s<\tau}\given\cF_s)+
                    \ind_{\tau\leq s}\bE(Y_t\ind_{\tau\leq s}\given\cG_s)\;.
\end{equation}
We then compute the two conditional expectations in (\ref{e:pro_expect}):
\begin{eqnarray*}
\bE(Y_t\ind_{s<\tau}\given\cF_s)
&=&\bE(Y_t\given\cF_s)-\bE(Y_t\ind_{\tau\leq s}\given\cF_s)\\
&=&\bE(m_t\given\cF_s)-\bE\big(\bE(\hat y_t(\tau)\ind_{\tau\leq s}\given\cF_t)\given\cF_s\big)\\
&=&m_s-\bE\big(\int_0^s\hat y_t(u)p_t(u)\nu(du)\given\cF_s\big)\\
&=&\tilde y_s G_s +\int_0^s\hat y_s(u)p_s(u)\nu(du)-\int_0^s\hat
y_s(u)p_s(u)\nu(du)=  \tilde y_s G_s \;,
\end{eqnarray*}   where we used Fubini-Tonelli's theorem and the condition \textit{(i)} to obtain the next-to-last identity.\\
Also, an application of Lemma \ref{l:pro_expect} yields to \
\begin{eqnarray*}
\bE(Y_t\ind_{\tau\leq s}\given\cG_s)
&=& \bE(\hat y_t(\tau)\ind_{\tau\leq s}\given\cG_s)
 =  \ind_{\tau\leq s}\frac{1}{p_s(\tau)}{\bE\big(\hat y_t(u)p_t(u)\given\cF_s\big)}_{\bgiven_{u=\tau}}\\
&=& \ind_{\tau\leq s}\frac{1}{p_s(\tau)}\hat y_s(\tau)p_s(\tau)=\ind_{\tau\leq s}\hat y_s(\tau)
\end{eqnarray*}
where the next-to-last identity holds in view of the condition
\textit{(ii)}. \finproof\brem \label{r:pro_CharMartLoc}
The extension of this characterization result to local martingales is more difficult.
Nevertheless, the sufficient condition holds: if $\big(\hat y_t(u)p_t(u),
t\geq u\big)$  and $(\tilde y_tG_t+\int_0^t\hat
y_t(u)p_t(u)\nu(du), t\geq 0)  $ are $(\bP,\bF)$-local
martingales, then using $\bF$-stopping times  $T_n$ for
localization, one gets that $(Y_{t\w T_n},t\geq 0)$ is a
$(\bP,\bG)$-martingale, hence $Y$ is a $(\bP,\bG)$-local
martingale.
\erem
\section{Canonical decomposition}\label{Sect:canonDecomp}

In this section, we work under $\bP$ and we  show that any
$\bF$-local martingale $x$ is a semimartingale in the initially
enlarged filtration $\bG^{\tau}$ and in the progressively enlarged
filtration $\bG$, and that   any  $ \bG$-local martingale is
a $ \bG^\tau$-semimartingale.  We also provide the canonical
decomposition of  any  $\bF$-local martingale as a semimartingale
in $\bG^\tau$ and in $\bG$. Under the assumption that the
$\bF$-conditional law of $\tau$ is absolutely continuous w.r.t.
the law of $\tau$, these questions were answered by Jacod in
\cite{Jacod85}, in the initial enlargement setting, and in
Jeanblanc and Le Cam \cite{JL_progressive}, in the
progressive enlargement case. Our aim here is to recover their
results in an alternative manner, under the $(\cE)$-Hypothesis.

We will need the following technical result, concerning the
existence of the predictable bracket $\langle x,p_.(u)\rangle$.
From Theorem 2.5.a in \cite{Jacod85}, it follows immediately that,
under the  $(\mathcal E)$-Hypothesis, for every
$(\bP,\bF)$-(local) martingale $x$, there exists a
$\nu$-negligible set $B$ (depending on $x$), such that $\langle
x,p_.(u)\rangle$ is well-defined for $u \notin B$. Hereafter,
by $\langle x,p_.(\tau)\rangle$ we mean ${\langle
x,p_.(u)\rangle}{\bgiven_{u=\tau}}$.

Furthermore, according to Theorem 2.5.b in \cite{Jacod85},
under the $(\mathcal E)$-Hypothesis, there exists an $\mathbb
F$-predictable increasing process $A$ and a $\mathcal P(\mathbb F)
\otimes \mathcal B(\mathbb R^+)$-measurable function $(t,\omega,u)
\rightarrow k_t(\omega,u)$ such that, for any $u \notin B$ and for
all $t \ge 0$,
\begin{equation}\label{Eq:crochetxpk}
{\langle x,p_.(u)\rangle}_t = \int_0^t k_s(u) p_{s-}(u) dA_s \quad \textrm{a.s.}
\end{equation}
(the two processes $A$ and $k$ depend on $x$, however, to simplify the notation, we do not write $A^{(x)}$, nor $k^{(x)}$).

Moreover,
\begin{equation}\label{Eq:integrk} \int_0^t \vert k_s(\tau) \vert
dA_s \; < \; \infty \quad \textrm{a.s., for any} \ t
>0.
\end{equation}
The following two propositions provide, under the $(\mathcal
E)$-Hypothesis, the canonical decomposition of any $(\bP,\bF)$-local martingale $x$ in the enlarged filtrations $\bG^{\tau}$ and
$\bG$, respectively. The case of initial enlargement has
been essentially established by Jacod (see Theorem 2.5.c in
\cite{Jacod85}), using a direct verification.
In our setting, one can obtain Jacod's result using the equivalent change of probability measure methodology
(see also Amendinger \cite{AmenTez}). Indeed, in view of Remark \ref{r:mart}, if $x$ is a
$(\bP,\bF)$-local martingale, it is a
$(\bP^*,\bG^\tau)$-local martingale, too. Noting that
$\frac{d\bP}{d\bP^*}=p_t(\tau)$ on $\cG^\tau_t$, Girsanov's
theorem tells us that the process $X^\tau$, defined by
$$
X^\tau_t := x_t-\int_0^t\frac{d\langle x,p_.(\tau)\rangle_s}{p_{s-}(\tau)}
$$
is a $(\bP,\bG^\tau)$-local martingale. However, the proof
presented here (for Proposition \ref{p:ini_Decompos}), is based on the $(\bP,\bG^\tau)$-martingales' characterization result given in Proposition \ref{p:ini_CharMart}.

\begin{pro}\label{p:ini_Decompos}\textbf{Canonical Decomposition in $\bG^{\tau}$}\\
Any $(\bP,\bF)$-local martingale $x$ is a $(\bP,\bG^\tau)$-semimartingale with canonical decomposition
$$
x_t=X^\tau_t+\int_0^t\frac{d\langle x,p_.(\tau)\rangle_s}{p_{s-}(\tau)},
$$
where $X^\tau$ is a $(\bP,\bG^\tau)$-local martingale.
\end{pro}
\proof  In view of Proposition \ref{p:ini_CharMart} and Remark \ref{r:ini_CharMartLoc}, using the notation $X^{\tau}=x(\tau)$, it suffices to show that, for $\nu$-almost every $u \ge 0$, the process
$$
x_t(u)p_t(u):=\Big(x_t-\int_0^t\frac{d\langle
x,p_.(u)\rangle_s}{p_{s-}(u)}\Big)p_t(u),\quad t\geq 0,
$$
is a $(\bP,\bF)$-local martingale. Indeed,
integration by parts formula leads to
$$
d(x_t(u)p_t(u))= p_{t-}(u) dx_t+x_{t-}(u) dp_t(u)+ d\big( \droit{x(u),p(u)}_t-\oblique{x(u),p(u)}_t\big).
$$
Hence, being the sum of three $(\bP,\bF)$-local martingales, the process $x(u)p(u)$ is a $(\bP,\bF)$-local martingale.
\finproof
Now, any $(\bP,\bF)$-local martingale is a $\bG$-adapted
process and a $(\bP,\bG^\tau)$-semimartingale (from the above
Proposition \ref{p:ini_Decompos}), so in view of Stricker's
theorem in \cite{Stricker}, it is also a $\bG$-semimartingale. The
following proposition aims to obtain the $\bG$-canonical
decomposition of an $\bF$-local martingale. We refer to
\cite{JL_progressive} for an alternative proof.

In order to study the canonical decomposition in $\bG$, we
add a regularity condition.
\begin{Ass}\label{a:regularite} There exists a version of the
process $(p_t(t),t\geq 0)$, such that $(\omega,t)\rightarrow
p_t(\omega,t)$ is $\cF_t\otimes\cB(\bR^+)$-measurable.
\end{Ass}
Then, the Az\'ema supermartingale $G$, introduced in Equation
(\ref{Eq:G}), admits   the Doob-Meyer decomposition  $G_t=\mu_t -
\int_0^t p_u(u)\nu(du), t \ge 0$, where $\mu$ is the $\mathbb
F$-martingale defined as
$$
\mu_t := 1-\int_0^t\left(p_t(u)-p_u(u)\right)\nu(du)
$$
(see, e.g., Section 4.2.1 in \cite{ejj}).

Before passing to the rigorous result on the canonical decomposition in $\bG$,
one can guess the form of the decomposition by means of a heuristic argument, based on
the equivalent change of probability measure: $(\bP,\bF)$-local martingale $x$  being a
$(\bP^*,\bG)$-local martingale,
$$x_t-\int_ 0^t \frac{1}{\ell_{s-}^*}d\oblique {x,\ell^*}_s$$ is a
$(\bP,\bG)$-local martingale, where $\ell^*:=\frac{1}{\ell} $ is
the Radon-Nikod\'ym density of $\bP$ w.r.t. $\bP^*$, given by
$$\frac{d\bP}{d\bP^*}\bgiven_{\cG_t}
=\ind_{t<\tau}\frac{G_t}{G^*(t)}+\ind_{\tau\leq
t}p_t(\tau)=\frac{1}{\ell_t}=\ell^*_t\;.$$
Based on the form of $\ell^*$, one has
$$
\ind_{s<\tau}d\oblique{x,\ell^*}_s=\ind_{s<\tau}\frac{d\oblique{x,G}_s}{G^*(s)}\;.
$$
This observation suggests Proposition \ref{p:pro_Decompos1} below, the proof of which is
based on the $(\bP,\bG)$-martingales' characterization result presented in Section \ref{Sect:martCaract}.
\begin{pro}\label{p:pro_Decompos1}\textbf{Canonical Decomposition in $\bG$}\\
Any (c\`{a}dl\`{a}g) $(\bP,\bF)$-local martingale $x$ is a
$(\bP,\bG)$-semimartingale with canonical decomposition
\begin{equation}\label{e:pro_Decompos1}
x_t= X_t+\int_0^{t\w\tau}\frac{d\langle
x,G\rangle_s}{G_{s-}}+\int_{t\w\tau}^t\frac{d\langle
x,p_.(\tau)\rangle_s}{p_{s-}(\tau)},
\end{equation} where $X$ is a
$(\bP,\bG)$-local martingale.
\end{pro}
\proof Relying on Remark \ref{r:pro_CharMartLoc}, we
check that $X$, defined in (\ref{e:pro_Decompos1}), is a
$(\bP,\bG)$-local martingale. We note that
$X_t=\ind_{t<\tau}\tilde x_t + \ind_{\tau\leq t}\hat x_t(\tau)$
with
$$
\tilde x_t = x_t-\int_0^t \frac{d\langle x,G\rangle_s}{G_{s-}}\;,\quad
\hat x_t(u)= x_t-\int_0^u \frac{d\langle x,G\rangle_s}{G_{s-}}
                -\int_u^t \frac{d\langle x,p_.(u)\rangle_s}{p_{s-}(u)}\;.
$$
We have to verify that
\begin{itemize}
\item [(i)] $(\hat x_t(u)p_t(u), t\geq u)$ is a $(\bP,\bF)$-local martingale;
\item [(ii)] $(\tilde x_tG_t + \int_0^t \hat x_t(u)p_t(u)\nu(du) , t\geq 0)$ is a $(\bP,\bF)$-local martingale.
\end{itemize}
In the proof of Proposition \ref{p:ini_Decompos}, we verified that (i) holds. In order to show (ii),
we apply It\^o's formula
\begin{eqnarray}
&& d\Big( \tilde x_tG_t + \int_0^t \hat x_t(u)p_t(u)\nu(du) \Big)\label{e:Ito_decom}\\
&& \quad\quad = G_{t-}dx_t + d(\droit{x,G}_t-\oblique{x,G}_t) + \tilde x_{t-}d\mu_t + (\hat x_t(t)-\tilde x_{t-})p_t(t)\nu(dt) + dz_t\;, \nonumber
\end{eqnarray}
where $z_t:=\int_{s=0}^t\int_{u=0}^s \nu(du) d\zeta_s(u)$ and
$\zeta_t(u):=\hat x_t(u)p_t(u)$. The first three terms on the
righthand-side of (\ref{e:Ito_decom}) are $(\bP,\bF)$-local
martingales, the fourth term is zero since $\hat x_t(t)- \tilde x_{t-}= \Delta \tilde x_t$,
and the $\bF$-adapted process $\tilde x$ has no jump at time $\tau$
(because, in our setting, $\tau$ avoids $\bF$-stopping times).
So, if we show that $z$ is a $(\bP,\bF)$-local martingale, we are done.
To do this, applying Fubini-Tonelli's theorem, one has
$$
z_t=\int_{u=0}^t\int_{s=u}^td\zeta_s(u)\nu(du)=\int_0^t (\zeta_t(u)-\zeta_u(u))\nu(du)\;.
$$
So, for $t_0<t$, one has
$$
z_t=\int_0^{t_0} (\zeta_t(u)-\zeta_u(u))\nu(du) + \int_{t_0}^t (\zeta_t(u)-\zeta_u(u))\nu(du)\;.
$$
Now, if $(\zeta_t(u),t\geq 0)$ is a $(\bP,\bF)$-martingale,
one gets $\bE(z_t\given\cF_{t_0})=z_{t_0}$. In the case where $\zeta_.(u)$ is a $(\bP,\bF)$-local martingale the
result is achieved by means of a localization argument.
\finproof 
\begin{rem} In a recent paper, Kchia et al.
\cite{Protter} have obtained the same decomposition formula, using
projection tools, in a more general setting. The main
challenge in their approach is that if a $\bG^\tau$-local
martingale is $\bG$-adapted, it is not necessarily a $\bG$-local
martingale (as remarked also by Stricker
\cite{Stricker}).\end{rem} The following lemma provides a formula
for the predictable quadratic covariation process
$\oblique{x,G}=\oblique{x,\mu}$ in terms of the density $p$.
\begin{pro}\label{l:crochet}
Let $x$ be a  $(\bP,\bF)$-local martingale and $\mu$ the $\mathbb
F$-martingale part in the Doob-Meyer decomposition of $G$. If
$kp_-$ is $dA \otimes d \nu$-integrable, then
\begin{equation}\label{e:crochet}
\oblique{x,\mu}_t=\int_0^t dA_s\int_s^\infty \nu(du)
k_s(u)p_{s-}(u) ,
\end{equation}
where $k$ was introduced in Equation (\ref{Eq:crochetxpk}).
\end{pro}
\proof First consider the right-hand side of (\ref{e:crochet}),
that is, by definition, predictable, and apply Fubini-Tonelli's theorem
\begin{eqnarray*}
&& \xi_t:= \int_0^t dA_s \int_s^\infty k_s(u)p_{s-}(u)\nu(du)\\
&& \quad = \int_0^t dA_s \int_s^t k_s(u)p_{s-}(u)\nu(du) + \int_0^t dA_s \int_t^\infty k_s(u)p_{s-}(u)\nu(du)\\
&& \quad = \int_0^t \nu(du) \int_0^u k_s(u)p_{s-}(u)dA_s + \int_t^\infty \nu(du) \int_0^t k_s(u)p_{s-}(u)dA_s\\
&& \quad = \int_0^t\oblique{x,p_{\cdot}(u)}_u \ \nu(du) + \int_t^\infty\oblique{x,p_{\cdot}(u)}_t \ \nu(du) \\
&& \quad = \int_0^\infty\oblique{x,p_{\cdot}(u)}_t \ \nu(du) +
          \int_0^t\left(\oblique{x,p_{\cdot}(u)}_u -\oblique{x,p_{\cdot}(u)}_t\right)\nu(du)\;.
\end{eqnarray*}
To verify (\ref{e:crochet}), it suffices to show that the process
$x\mu-\xi$ is an $\bF$-local martingale (since $\xi$ is a
predictable, finite variation process). By definition, for
$\nu$-almost every $u\in\mathbb R^+$, the process
$\left(m_t(u):=x_t p_t(u)-\oblique{x,p_{\cdot}(u)}_t, t\geq
0\right)$ is an $\bF$-local martingale. Then, given that
$1=\int_0^{\infty} p_t(u)\nu(du)$ for every $t \ge 0$, a.s., we
have
\begin{eqnarray*}
&& x_t\mu_t-\xi_t = x_t\int_0^\infty p_t(u)\nu(du) - x_t\int_0^t\left(p_t(u)-p_u(u)\right)\nu(du)\\
&& \quad  -\int_0^\infty\oblique{x,p_{\cdot}(u)}_t \ \nu(du) +
            \int_0^t\left(\oblique{x,p_{\cdot}(u)}_t-\oblique{x,p_{\cdot}(u)}_u\right)\nu(du)\\
&& \quad = \int_0^\infty m_t(u)\nu(du) - \int_0^t
\left(m_t(u)-m_u(u)\right)\nu(du) + x_t\int_0^t p_u(u) \nu(du)-
\int_0^t p_u(u) x_u \nu(du)\;.
\end{eqnarray*}
The first two terms are local martingales, in view of the
martingale property of $m(u)$. As for the last term, using the
fact that $\nu$ has no atoms, we find
\begin{eqnarray*}
\lefteqn{d\left( x_t \int_0^t p_u(u)\nu(du) - \int_0^t p_u(u)
x_u\nu(du)\right)} \\&=& \left( \int_0^t p_u(u)\nu(du)
\right)\,dx_t+ x_t p_t(t)\nu(dt) - p_t(t) x_t\nu(dt) =
\left(\int_0^t p_u(u)\nu(du)\right) \ dx_t
\end{eqnarray*}
and we have, indeed, proved that $x\mu -\xi$ is an $\mathbb
F$-local martingale. 
\finproof
We end this section proving that any
$(\bP^*,\bG)$-martingale remains a
$(\bP^*,\bG^{\tau})$-semimar\-tingale, but it is not necessarily a
$(\bP^*,\bG^\tau)$-martingale. Indeed, we have the following
result.
\begin{pro}\label{23}
Any $(\bP^*,\bG)$-martingale  $Y^*$ is a $(\bP^*,\bG^\tau)$-
semimartingale which can have a  non-null bounded variation part.
\end{pro}
\proof The result follows immediately from Proposition
\ref{p:pro_CharMart} (under $\mathbb P^*$), noticing that the
$(\bP^*,\bG)$-martingale $Y^*$ can be written as $Y_t^*=\tilde
y_t^* \ind_{t<\tau}+\hat y_t^*(\tau)\ind_{\tau\leq t}$. Therefore,
in the filtration $\mathbb G^{\tau}$, it is the sum of two
$\mathbb G^{\tau}$-semimartingales: the processes  $\ind_{t<\tau}
$ and $\ind_{\tau\leq t}$ are 
$\mathbb G^{\tau}$-semimartingales, as well as the processes
$\tilde y, \hat y ^*(\tau)$.  Indeed, from Proposition
\ref{p:pro_CharMart},   recalling that the $(\mathbb P^*,\mathbb
F)$-density of $\tau$ is a constant equal to one, we know that,
for every $u>0 $, $\big(\hat y_t^*(u), t\geq u\big)$ is an
$\mathbb F$-martingale and   that the process $\big(\tilde y_t^*
G^*(t)+\int_0^t\hat y_u^*(u)\nu(du),t \geq 0 \big)$ is an $\mathbb
F$-martingale, hence $\tilde y^*$ is a $\mathbb
G$-semimartingale.

It can be noticed that the  $(\bP^*,\mathbb G)$-martingale $M^*$,
defined in (\ref{Eq:M*}), is such that $M^*_t$ is, for any $t$,
a $\mathcal G^\tau_0$-measurable random variable. Therefore, $M^*
$ is not a $(\bP^*,\mathbb G^\tau )$-martingale, since, for $s \le
t$, $\bE( M^*_t \vert \mathcal G_s^\tau)= M^*_t \neq  M^*_s $, but
it is a bounded variation $\bG^\tau$-predictable process, hence a
$\mathbb G^\tau$-semimartingale with null martingale part.  In
other terms, $\mathbb H$ is not immersed in $\mathbb  G^\tau$
under $\bP^*$. \finproof

As in  Lemma  \ref{23}, we deduce that  any $(\bP,\bG)$-martingale
is a $(\bP,\bG^\tau)$-semimartingale. Note that this result can
also be proved using Lemma \ref{23} and a change of probability
argument: a $(\bP,\bG)$-martingale is a
$(\bP^*,\bG)$-semimartingale (from Girsanov's theorem), thus also
a $(\bP^*,\bG^\tau)$-semimartingale in view of Lemma \ref{23}. By
another use of Girsanov's theorem, it is thus a
$(\bP,\bG^\tau)$-semimartingale.

\section{Predictable Representation Theorems}\label{Sect:PRP}
The aim of this section is to obtain Predictable Representation
Property (PRP   hereafter) in the enlarged filtrations $\mathbb G$
and $\mathbb G^{\tau}$, both under $\mathbb P$ and $\mathbb P^*$.
To this end, we assume that there exists a $(\bP,\bF)$-local
martingale $z$ (possibly multidimensional), such that the PRP holds in
$(\bP,\bF)$ (cf. Assumption \ref{a:PRP_F}, below). Notice that $z$ is not necessarily continuous.

Beforehand we introduce some notation: ${\mathcal
M_{\mathrm{loc}}}(\bP,\bF)$ denotes the set of $(\bP,\bF)$-local
martingales, and ${\mathcal M}^2(\bP,\bF)$ denotes the set of
$(\mathbb P,\mathbb F)$-martingales $x$, such that
\begin{equation}\label{e:int_condition}
\mathbb E \left( x_t^2 \right) < \infty, \quad \forall \ t \ge 0.
\end{equation}

Also, for a $(\bP,\bF)$-local martingale
$m$, we denote by $\cL(m,\bP,\bF)$ the set of $\bF$-predictable
processes which are integrable with respect to $m$ (in the sense
of local martingale), namely (see, e.g., Definition 9.1 and
Theorem 9.2. in He \emph{et al.} \cite{HeWangYan})
$$
\cL(m,\bP,\bF)=\left\{\varphi\in\cP(\bF):\ \left(\int_0^{\cdot} \varphi^2_sd[m]_s \right)^{1/2}\ \mbox{is}\ \bP-\mbox{locally integrable} \right\}\;.
$$

\begin{Ass}\label{a:PRP_F}\textbf{PRP for $(\bP,\bF)$}\\
There exists a process $z\in {\mathcal M_{\mathrm{loc}}}(\bP,\bF)$
such that every $x\in {\mathcal M_{\mathrm{loc}}}(\bP,\bF)$ can be
represented as
$$x_t=x_0+\int_0^t\varphi_s dz_s$$
for some $\varphi\in\cL(z,\bP,\bF)$.
\end{Ass}

We start investigating what happens under the measure $\bP^*$, in
the initially enlarged filtration $\mathbb G^{\tau}$.

Recall that, assuming the immersion property, Kusuoka
\cite{Kusuoka} has established a PRP for the progressively
enlarged filtration, in the case where $\bF$ is a Brownian
filtration.

Also, under the \textit{equivalence assumption} in $[0,T]$ and
assuming a PRP in the reference filtration $\mathbb F$, Amendinger
(see \cite[Thm. 2.4]{AmenTez}) proved a PRP in $(\mathbb
P^*,\mathbb G^{\tau})$ and extended the result to $(\mathbb
P,\mathbb G^{\tau})$, in the case where  the underlying (local)
martingale in the reference filtration is continuous.
Under the $(\cE)$-Hypothesis,  Grorud and Pontier \cite[Prop. 4.3]{GrPontAsym} have established
a PRP for $(\bP,\bG^\tau)$-local martingales, in the case where the filtration $\bF$ consists
of a point process and a continuous martingale (typically a Brownian motion).

\begin{pro}\label{p:ini_PRP}\textbf{PRP for $(\bP^*,\bG^\tau)$}\\
Under Assumption \ref{a:PRP_F}, every $X^\tau\in {\mathcal
M}_{\mathrm{loc}}(\bP^*,\bG^\tau)$ admits a representation
\begin{equation}\label{e:prt_Xtau}
X^\tau_t=X^\tau_0+\int_0^t\Phi^\tau_s dz_s
\end{equation}
where $\Phi^{\tau}\in \cL(z,\bP^*,\bG^\tau) $. In the case where
$X^\tau\in {\mathcal M}^2(\bP^*,\bG^\tau)$, one has $\bE^*\big(\int_0^t {(\Phi_s^{\tau})}^2 d[
z]_s\big)<\infty$, for all $t \ge 0$ and the representation is
unique.
\end{pro}
\proof From Theorem 13.4 in \cite{HeWangYan}, it suffices to prove
that any bounded   martingale admits a predictable representation
in terms of $z$. Let $X^\tau\in {\mathcal
M_{\mathrm{loc}}}(\bP^*,\bG^\tau)$ be bounded by $K$. From
Proposition \ref{p:ini_CharMart},  $X^\tau_t=x_t(\tau)$  where,
for $\nu$-almost every $u\in\mathbb R^+$, the process
$\big(x_t(u),t\geq 0\big)$ is a $(\bP^*,\bF)$-martingale, hence a
$(\bP,\bF)$-martingale. Thus Assumption \ref{a:PRP_F} implies that
(for $\nu$-almost every $u\in\mathbb R^+$),
$$
x_t(u)=x_0(u)+\int_0^t\varphi_s(u) dz_s\;,
$$
where $(\varphi_t(u),t\geq 0)$ is an $\bF$-predictable process.

The process $X^\tau$ being bounded by $K$, it follows by an application of Lemma
\ref{l:ini_expect}(i) that for $\nu$-almost every $u\geq 0$, the process $(x_t(u),t\geq 0)$ is bounded by $K$.
Then, using the It\^o isometry,
\begin{eqnarray*}
\bE^*\big( \int_0^t \varphi^2_s(u)d\droit{z}_s\big)
&=&\bE^*\big(\int_0^t\varphi_s(u)dz_s\big)^2\\
&=&\bE^*\big((x_t(u)-x_0(u))^2\big)\leq \bE^*(x^2_t(u) )\leq K^2\;.
\end{eqnarray*}
Also, from Stricker and Yor \cite[Lemma 2]{StrickerYor}, one
can consider a version of the process $\int_0^\cdot \varphi^2_s(u)
d[z]_s$ which is measurable with respect to $u$. Using this fact,
$$
\bE^*\Big[\big(\int_0^t \varphi _s^2(\tau) d [z]_s\big)^{1/2}\Big]=
\int_0^\infty \nu(du) \Big(\bE^*\big(\int_0^t \varphi _s^2(u) d [z]_s\big)\Big)^{1/2}\leq
\int_0^\infty \nu(du)   K =K\;.
$$

The process $\Phi^\tau$ defined by $\Phi^\tau_t=\varphi_t(\tau)$
is $\bG^\tau$-predictable, according to Proposition
\ref{p:ini_categorize}, it satisfies (\ref{e:prt_Xtau}), with
$X_0(\tau)=x_0(\tau)$, and it belongs to $\cL(z,\bP^*,\bG^\tau)$.

If $X^\tau\in {\mathcal M}^2(\bP^*,\bG^\tau)$, from It\^o's isometry,
$$
\bE^*\left(\int_0^t(\Phi^\tau_s)^2 d[z]_s\right)=
\bE^*\left(\int_0^t\Phi^\tau_s dz_s\right)^2=\bE^* (X^\tau_t -
X_0^{\tau})^2 <\infty\;.
$$
Also, from this last equation, if $X^\tau\equiv 0$ then
$\Phi^\tau\equiv 0$, from which the uniqueness of the
representation follows. \finproof

Passing to the progressively enlarged filtration $\bG$, which
consists of two filtrations, $\mathbb G=\mathbb F \vee \mathbb H$,
intuitively one needs two martingales to establish a PRP. Apart
from $z$, intuition tells us that a candidate for the second
martingale might be the compensated martingale  of $H$, that was
introduced, respectively under $\mathbb P$ (it was denoted by $M$)
and under $\mathbb P^*$ (denoted  by  $M^*$), in Equation
(\ref{Eq:M}) and in Equation (\ref{Eq:M*}).

\begin{pro}\label{p:pro_PRP}\textbf{PRP for $(\bP^*,\bG)$}\\
Under Assumption \ref{a:PRP_F}, every $X\in {\mathcal
M}_{\mathrm{loc}}(\bP^*,\bG)$ admits a representation
$$
X_t=X_0+\int_0^t \Phi_s dz_s + \int_0^t \Psi_s dM^*_s
$$
for some processes $\Phi\in\cL(z,\bP^*,\bG)$ and
$\Psi\in\cL(M^*,\bP^*,\bG)$. Moreover, if $X\in {\mathcal M}^2(\bP^*,\bG)$,
 one has, for any $t \ge0$,
$$
\bE^* \left( \int_0^t \Phi_s^2 d[ z]_s \right) < \infty\quad ,
\quad \bE^* \left( \int_0^t  \Psi_s^2 \lambda^*(s) \nu(ds) \right)
< \infty\;,
$$
and the representation is unique.
\end{pro}
\proof It is known that any $(\bP^*,\bH)$-local martingale $\xi$
can be represented as $\xi_t=\xi_0+\int_0^t\psi_s dM^*_s$ for some
process $\psi\in\cL(M^*,\bP^*,\bH)$
(see, e.g., the proof in Chou and Meyer \cite{ChouMeyer}). Notice
that $\psi$ has a role only before $\tau$ and, for this reason, $\psi$ can be chosen deterministic.

Under $\bP^*$, we then have\vspace{-1mm}
\begin{itemize}
\item the PRP holds in $\bF$ with respect to $z$,\vspace{-1mm}
\item the PRP holds in $\bH$ with respect to $M^*$,\vspace{-1mm}
\item the filtration $\bF$ and $\bH$ are independent.\vspace{-1mm}
\end{itemize}
From classical literature (see  Lemma 9.5.4.1(ii) in
Jeanblanc \emph{et al.} \cite{ChJeYor09}, for instance)
the filtration $\bG=\bF\V\bH$ enjoys the PRP under $\bP^*$ with
respect to the pair $(z,M^*)$.

Now suppose that $X\in {\mathcal M}^2(\bP^*,\bG)$. We find
\begin{eqnarray*}
\infty &>&\bE ^*{(X_t - X_0)}^2
= \bE^*\left(\int_0^t \Phi_s dz_s + \int_0^t \Psi_s dM^*_s\right)^2\\
&=& \bE^* \left(\int_0^t\Phi^2_s d[z]_s \right)+ 2\bE^*
\left(\int_0^t \Phi_s dz_s \int_0^t \Psi_s
dM^*_s\right)+\bE^*\left(\int_0^t\Psi^2_s
\lambda^*(s)\nu(ds)\right),
\end{eqnarray*}
where in the last equality we used the It\^o isometry. The
cross-product term in the last equality is zero due to the
orthogonality of $z$ and $M^*$ (under $\bP^*$). From this
inequality, the desired integrability conditions hold and the
uniqueness of the representation follows (as in the previous
proposition). \finproof

\begin{rem}
In order to establish a PRP for  the  initially enlarged
filtration $\bG^\tau$ and under $\bP^*$, one could have proceeded
as in the proof of Proposition \ref{p:pro_PRP}, noting that any
martingale $\xi   $ in the ``constant'' filtration $\sigma(\tau)$
satisfies
 $\xi_t=\xi_0+0$ and that under $\bP^*$ the two
filtrations $\bF$ and $\sigma(\tau)$ are independent.
\end{rem}

\begin{pro}\label{p:PRP}\textbf{PRP under $\bP$}\\ Under Assumption \ref{a:PRP_F}, one has:
\begin{itemize}
\item [(i)] Every $X^\tau\in {\mathcal
M_{\mathrm{loc}}}(\bP,\bG^\tau)$ can be represented as
$$
X^\tau_t=X_0^{\tau}+\int_0^t \Phi^\tau_s dZ^{\tau}_s\;
$$
where $Z^\tau$ is the martingale part in the $\bG^\tau$-canonical decomposition of $z$
and $\Phi\in\cL(Z^\tau,\bP,\bG^\tau)$.
\item[(ii)] Every $X\in {\mathcal M_{\mathrm{loc}}}(\bP,\bG)$ can
be represented as
$$
X_t=X_0+\int_0^t \Phi_s dZ_s + \int_0^t \Psi_s dM_s,
$$
where $Z$ is the martingale part in the $\bG$-canonical decomposition
of $z$ (cf. Equation \ref{e:pro_Decompos1}), $M$ is the $(\bP,\bG)$-compensated martingale associated with
$H$ and $\Phi\in\cL(Z,\bP,\bG)$, $\Psi\in\cL(M,\bP,\bG)$.
\end{itemize}
\end{pro}
\proof The assertion (i) (resp. (ii)) follows from Proposition
\ref{p:ini_PRP} (resp. Proposition \ref{p:pro_PRP}) and the
stability of PRP  under an equivalent change of measure (see for
example Theorem 13.12 in \cite{HeWangYan}).
\finproof
The PRP for the progressively enlarged filtration with a random
time (i.e., part (ii) of the above proposition), has been first
presented by Jeanblanc and Le Cam \cite{JL_immersion}. Our proof
has the advantage that it can be straightforwardly generalized to
the case where $\tau$ is a vector of random times, as is discussed
in the last section.
\section{Concluding Remarks}\label{Sect:conclusion}
We conclude the paper with some important comments:
\begin{itemize}
\item In the multidimensional case, that is when $\tau=(\tau_1,\dots,\tau_d)$ is a vector of finite
random times, the same machinery can be applied. More precisely, under the assumption
$$
\bP(\tau_1\in d\theta_1, \dots , \tau_d\in\ d\theta_d\given\cF_t)
\sim \bP(\tau_1\in d\theta_1, \dots , \tau_d\in d\theta_d)
$$
one defines the probability $\bP^*$ equivalent to  $\mathbb P$ on
$\cG^\tau_t=\cF_t\V\sigma(\tau_1)\V\dots\V\sigma(\tau_d)$ by
\begin{equation}\label{e:MultiDim_P*}
{\frac{d\bP^*}{d\bP}}\bgiven_{\mathcal G_t^{\tau}}=\frac{1}{p_t(\tau_1,\dots,\tau_d)},
\end{equation}
where $p_t(\tau_1,\dots,\tau_d)$ is the (multidimensional) analog
to $p_t(\tau)$, and the results for the initially enlarged
filtration are obtained in the same way as for the one-dimensional case.

As for the progressively enlarged filtration, we define $\bH:=\bH^1\V\cdots\V\bH^d$, where
$\bH^i$ stands for the natural filtration of the indicator process $H^i=(\ind_{\tau_i\leq t},t\geq 0)$.
The progressive enlargement of $\bF$ with the vector $(\tau_1,\dots,\tau_d)$, is then defined by
(the right-continuous regularization of) the filtration $\bF\V\bH$.
One has to note that, in this case, a measurable process is decomposed into $2^d$ terms,
corresponding to the measurability of the process on the various sets
$\{\tau _i \leq t <\tau_j, i\in I, j\in I^c\}$ for all the subsets $I$ of $ \{1,...,d\}$.

An interesting point is the generalization of the proof of predictable
representation theorem for the progressively enlarged filtration, in the multidimensional case.
Under the probability $\bP^*$, defined in (\ref{e:MultiDim_P*}), the filtration $\bH$ is
independent of $\bF$. So, once a PRP for $\bH$ holds true, it is
straightforward to generalize Proposition \ref{p:pro_PRP} to the multidimensional case.
To this end, we further assume that the (joint) law
$\bP(\tau_1\in d\theta_1, \dots , \tau_d\in d\theta_d)$ of $\tau$ is absolutely continuous w.r.t. the Lebesgue measure on $(\bR^+)^d$.
One then has $\bP(\tau_i=\tau_j)=0$ and thus the process
$H=(H^1,...,H^n)$ is an $n$-variate point process (in the terminology of Br\'emaud \cite{Bremaud}).
So, $\bH$ enjoys the PRP with respect to the compensated martingales of
$H^1,...,H^n$ (see for instance Br\'emaud \cite{Bremaud} Chap III, Sec 3, Theorems 9 \& 11).

\item In this study, \emph{honest} times (recall that a random
time $L$ is honest if, for any $t$, it is equal to an $\mathcal
F_t$-measurable random variable on $\{ L < t \}$) are
automatically excluded, as we explain now. Under the probability
$\bP^*$, the Az\'ema supermartingale associated with $\tau$, being
a continuous decreasing function, has a trivial Doob-Meyer
decomposition $G^*=1-A^*$ with $A^*_t=\int_0^t\nu(du)$. So,
$A^*_\infty=1$ and, in particular, $\tau$ can not be an honest
time: recall that in our setting, $\tau$ avoids the $\bF$-stopping
times and therefore, from a result due to Az\'ema (cf.
\cite{Azema}, part (b) of the Theorem on pages 300 and 301), if
$\tau$ is an honest time, the random variable $A^*_\infty$ should
have exponential law with parameter 1, which is not the case (
note that the notion of honest time does not depend on the
probability measure). \item  Under the $(\mathcal E)$-Hypothesis,
the immersion property between $\bF$ and $\bG$ is equivalent to
$p_t(u)=p_u(u), t\geq u$ (cf. \cite[Corollary 3.1]{JL_progressive}). In particular, as expected, the canonical
decomposition formula presented in Proposition
\ref{p:pro_Decompos1} is trivial, that is, the two integral terms
on the right-hand side of (\ref{e:pro_Decompos1}) vanish. \item
Predictable representation  theorems can be obtained in the more
general case, where any $(\bP,\mathbb F)$-martingale $x$ admits a
representation as $$x_t = x_0+\int _0^t \int _E \varphi (s,\theta)
\tilde \mu (ds,d\theta),$$ for a compensated martingale associated
with a point process.
\end{itemize}

\begin{center}
{\bf Acknowledgements}
\end{center}
The authors would like   to thank Prof. Sh. Song (UEVE, France) and Prof. M. Pontier
(Paul Sabatier, Toulouse, France) for their enlightening discussions and precious remarks.

\end{document}